\documentclass[12pt]{article}
\usepackage{epsfig}
\usepackage{algorithm}
\usepackage{algorithmic}
\usepackage{latexsym}
\usepackage{amsmath}
\usepackage{cite}
\usepackage{verbatim}
\usepackage{amsfonts}
\usepackage{graphicx}
\usepackage[normalem]{ulem}
\usepackage{array}
\usepackage{enumerate}
\usepackage{authblk}
\usepackage{color}\usepackage[dvipsnames]{xcolor}
\usepackage{multirow}
\usepackage{tikz}\usetikzlibrary{patterns}
\usepackage{makecell, amssymb, microtype}
\usepackage{theorem}

\DeclareFontFamily{OT1}{pzc}{}
\DeclareFontShape{OT1}{pzc}{m}{it}{<-> s * [1.10] pzcmi7t}{}
\DeclareMathAlphabet{\mathpzc}{OT1}{pzc}{m}{it}

\usepackage[
  bookmarks=false, 
  colorlinks,
  citecolor=brown!70!black,
  linkcolor=brown!80!black,
  urlcolor=blue!70!black,
]{hyperref}
\usepackage{relsize,exscale}
\newtheorem{lem}{Lemma}[section]
\newtheorem{thm}{Theorem}[section]
\newtheorem{prop}{Proposition}[section]

\theorembodyfont{\rm}
\newtheorem{defn}{Definition}[section]
\newtheorem{rem}{Remark}[section]

\newcommand\oeis[1]{\href{https://oeis.org/#1}{#1}}

\providecommand{\keywords}[1]{\textbf{Keywords:} #1}

\let\le\leqslant

\let\leq\leqslant
\let\geq\geqslant

\date{\today}
\author[1]{Jean-Luc Baril}
\author[1]{Sergey Kirgizov}
\author[2]{Armen Petrossian}
\affil[1]{\rm LIB, Université de Bourgogne Franche-Comté\protect\\
  B.P. 47 870, 21078 Dijon Cedex France\protect\\
   {\tt E-mail: \{barjl,sergey.kirgizov\}@u-bourgogne.fr
   }
}
\affil[2]{\rm ESEO Paris\protect\\
78140 Vélizy-Villacoublay\protect\\
   {\tt E-mails: armen.petrossian@eseo.fr
   }
}
\date{\today}

\title{Dyck paths with catastrophes modulo the positions of a given pattern}

\begin{document}

\maketitle
\begin{abstract}
For any pattern $p$ of length at most 2, we provide generating
functions and asymptotic approximations for the number of
$p$-equivalence classes of Dyck paths with catastrophes, where two
paths of the same length are $p$-equivalent whenever the positions of
the occurrences of the pattern $p$ are the same.
\end{abstract}
\section{Introduction and notation}

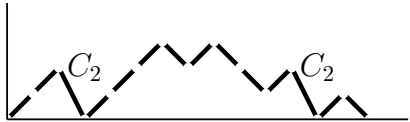
\begin{figure}[h]
\begin{center}
\scalebox{0.86}{\begin{tikzpicture}[ultra thick]
   \draw[black, line width=2pt](0,0)--(0.4,0.4)--(0.8,0.8)--(1.2,0)--(1.6,0.4)--(2,0.8)--(2.4,1.2)--(2.8,0.8)--(3.2,1.2)--(3.6,0.8)--(4,0.4)--(4.4,0.8)--(4.8,0)--(5.2,0.4)--(5.6,0);
   \draw  (1.2,0.8) node {\large $C_2$};\draw  (4.8,0.8) node {\large $C_2$};

   \tikzset{every node/.style={circle, white, fill = white, inner sep = 1.5pt}}
   \node at (0,0) {};
   \node at (0.4,0.4) {};
   \node at (0.8,0.8) {};
   \node at (1.2,0) {};
   \node at (1.6,0.4) {};
   \node at (2,0.8) {};
   \node at (2.4,1.2) {};
   \node at (2.8,0.8) {};
   \node at (3.2,1.2) {};
   \node at (3.6,0.8) {};
   \node at (4,0.4) {};
   \node at (4.4,0.8) {};
   \node at (4.8,0) {};
    \node at (5.2,0.4) {};
    \node at (5.6,0) {};
   \draw[black, thick] (0,0)--(6.2,0);\draw[black, thick] (0,0)--(0,1.8);
\end{tikzpicture}}
\end{center}
\caption{Dyck path with catastrophes $UUC_2UUUDUDDUC_2UD$. }
\label{fig1}
\end{figure}

\noindent A \emph{Dyck path with catastrophes} is a lattice path in
the first quadrant of the $xy$-plane that starts at the origin, ends
on the $x$-axis, and is made of up-steps $U=(1,1)$, down-steps
$D=(1,-1)$, and catastrophe steps $C_k=(1,-k)$, $k\geq 2$, so that
catastrophe steps always end on the $x$-axis. Depending on the
context, we can use the symbol $C$ to design a catastrophe step, and
by convenience we use $C_1=D$.  We let $\mathcal{E}$ denote the set of
all Dyck paths with catastrophes, and $\mathcal{D}$ be the set of
\emph{Dyck paths}, i.e.~the paths in $\mathcal{E}$ that do not contain
any catastrophe steps $C_k$, $k\geq 2$.  The \emph{length} $|P|$ of a
path $P$ is the number of its steps. The empty path is denoted by
$\epsilon$.  See Figure \ref{fig1} for an example of a Dyck path with
catastrophes of length $14$. A {\it pattern} consists of consecutive
steps in a path. We will say that {\it an occurrence of a pattern} (or
for short {\it a pattern}) is at {\it position} $i\geq 1$ in a path
whenever the first step of the pattern appears at the $i$-th step of
the path, the second step at the $(i+1)$-st step, and so on. The {\it
  height} of an occurrence of a pattern is the minimal ordinate
reached by its points. For instance, the Dyck path with catastrophes
$P=UUC_2UUUDUDDUC_2UD$ contains three occurrences of the pattern $UU$
at positions $1$, $4$ and $5$, and the heights of these occurrences
are respectively $0$, $0$ and $1$.

The concept of a Dyck path with catastrophes was first introduced by
Krinik et al. in~\cite{Kri} in the context of queuing theory. These
paths correspond to the evolution of a queue by allowing some resets
modeled by a catastrophe step $C_k$, $k\geq 2$. Then, Banderier and
Wallner~\cite{Ban} provided enumerative results and limit laws of these
objects. They showed how any non empty path $P\in\mathcal{E}$ can be
decomposed either as $P=U\alpha D\beta$, or
$P=U\alpha_1U\alpha_2\ldots U\alpha_kC_k\beta$ for some $k\geq 2$,
where $\alpha, \alpha_1,\alpha_2, \ldots,\alpha_k$ are Dyck paths in
$\mathcal{D}$ and $\beta\in\mathcal{E}$. They deduced a functional
equation for the generating function $E(x)=\sum_{n\geq 0} e_n x^n$
where $e_n$ is the number of paths of length $n$ in $\mathcal{E}$,
with the solution $$E(x)=\frac {2x-1+\sqrt
  {1-4\,{x}^{2}}}{x-1+(1+x)\sqrt {1-4\,{x}^{2}}}.
$$ The sequence $(e_n)_{n\geq 0}$ corresponds to \oeis{A224747} in the
On-line Encyclopedia of Integer Sequences (OEIS)~\cite{OEIS}, and the
first values for $n\geq 0$ are $1,0,1,3,5,12,23,52,105,232,480$.
More
recently, Baril and Kirgizov~\cite{Bari4} exhibited a one-to-one
correspondence between Dyck paths with catastrophes of length $n$ and
Dyck paths of length $2n$ avoiding $UUU$ and $DUD$ at height at least
one, and where every occurrence of $UD$ on the $x$-axis appears
before (but not necessarily contiguous with) an occurrence of $UUU$.

On the other hand, in~\cite{Bari1,Bari2,Bari4,Bari3,Barasii,Manes} the authors
investigated equivalence relations on the sets of Dyck paths, Motzkin
paths, skew Dyck paths, \L{}ukasiewicz paths, and Ballot paths where
two paths of the same length are equivalent whenever they coincide on
all occurrences of a given pattern. The main goal of this study
consists in extending these studies for Dyck paths with catastrophes
by considering the analogous equivalence relation on $\mathcal{E}$.

\begin{defn}
Two Dyck paths of the same length and with catastrophes are
$p$-equivalent whenever they have the same positions of the
occurrences of the pattern $p$.
\end{defn}
For instance, the path $UUDUUC_3$ is $U$-equivalent to $UUC_2UUC_2$
since the occurrences of $U$ appear at the same positions in the two
paths.

Then, one may naturally split the set $\mathcal{E}$ into {\it
  $p$-equivalence classes}, which are constructed so that paths $P$
and $Q$ in $\mathcal{E}$ belong to the same $p$-equivalence class if,
and only if, they are $p$-equivalent.

In this paper, we provide ordinary generating functions (o.g.f.\ for
short) for the number of $p$-equivalence classes in $\mathcal{E}$ with
respect to the length whenever $p$ is a pattern of length at most two.
Our method consists in providing one-to-one correspondences between
equivalence classes and certain subsets of $\mathcal{E}$ (called
subsets of {\em representative elements}) and enumerating them using
algebraic techniques.  Remark that only one pattern of size 2, namely
$DD$, gives a non-rational generating function. For this pattern, the
construction of a set of representative elements and its enumeration
are quite intricate and handled in the last subsection of the paper.
We refer to Table~\ref{tt} for an overview of our results.

\begin{table}[H]
\begin{center}
\resizebox{13.8cm}{!}{\begin{tabular}{| c | c | l| c|c|}
\hline
Pattern  & OEIS~\cite{OEIS} & Sequences for $a_n$, $2 \leq n \leq 10$, and formul\ae & Theorem& O.g.f.\ \\ \hline
$U$ & \mbox{Shift of }\oeis{A037952} &
 \makecell[l]{
   1, 1, 3, 4, 10, 15, 35, 56, 126, ...\\
 $\displaystyle {n-1 \choose {\lfloor\frac{n-2}{2}\rfloor}}$
 }
&Thm~\ref{thm2} & \multirow{2}{*}{Algebraic}
 \\ \cline{1-4}
 $D$  & New   &
 \makecell[l]{
   1, 1, 3, 5, 11, 21, 42, 84, 162, ... \\
   $\begin{cases}
     \displaystyle
         {n \choose \frac{n-3}{2}} &\text{ if } n \text{ is odd,} \\
         \displaystyle
             {n \choose  \frac{n-4}{2} }+
             {n \choose \frac{n}{2}}\slash\left(\frac{n}{2}+1\right) &\text{ otherwise}
 \end{cases}$
 }
&Thm~\ref{thm4} &  \\ \hline
 $C$ &\oeis{A212804}  &
 \makecell[l]{
   1, 1, 2, 3, 5, 8, 13, 21, 34, ... \\
   $\begin{cases} c_{n}=c_{n-1}+c_{n-2},\\ c_0=1, c_1=0\end{cases}$
 }
& Thm~\ref{thm6} & \multirow{7}{*}{}
  \\ \cline{1-4}
  $UU$ & \oeis{A347493}  &
  \makecell[l]{
    1, 1, 3, 4, 8, 13, 24, 41, 75, ... \\  
    $\begin{cases}f_{n}=f_{n-1}+f_{n-2}+f_{n-4},\\  f_0=1, f_1=0, f_2=1, f_3=1\end{cases}$
      }
&Thm~\ref{thm8}  &
   \\ \cline{1-4}
   $UD$ &  \oeis{A215004} &
   \makecell[l]{
      1, 1, 3, 5, 10, 17, 30, 50, 84, ... \\
      $\begin{cases}g_n = 2g_{n-1} + g_{n-2} - 3g_{n-3} + g_{n-5},\\g_1=0, g_0=g_2=g_3=1, g_4=3\end{cases}$
   }
&Thm~\ref{thm10}&
    \\ \cline{1-4}
    $UC$ &  New &
    \makecell[l]{
       1, 1, 2, 4, 5, 9, 15, 24, 40, ...\\
       $\begin{cases}i_n =2i_{n-1}-i_{n-2}+i_{n-3}-i_{n-4}+i_{n-5}-i_{n-6}+i_{n-7},\\ i_0=1, i_1=0, i_2=1, i_3=1, i_4=2, i_5=4, i_6=5\end{cases}$
    }
 &Thm~\ref{thm12}   &Rational
 \\  \cline{1-4}
    $DC$ &  New &
    \makecell[l]{
      1, 1, 1, 2, 2, 4, 6, 9, 14, ...\\
      $\begin{cases}j_n=2j_{n-1}-j_{n-2}+j_{n-5}-j_{n-6}+j_{n-7},\\ j_0=1, j_1=0, j_2=j_3=j_4=1, j_5=2, j_6=2\end{cases}$
    }
   & Thm~\ref{thm14}&
    \\ \cline{1-4}
    $CU$ & \mbox{Shift of }  \oeis{A000045} &
    \makecell[l]{
      1, 1, 1, 2, 3, 5, 8, 13, 21, ...\\
      $\begin{cases}k_n=k_{n-1}+k_{n-2},\\ k_0=1, k_1=0, k_2=k_3=k_4=1\end{cases}$
    }
&  Thm~\ref{thm16}  &  \\
\cline{1-4}
$DU$ & \oeis{A212804}   &
\makecell[l]{
  1, 1, 2, 3, 5, 8, 13, 21, 34, ...\\
  $\begin{cases}\ell_n = \ell_{n-1} + \ell_{n-2},\\\ell_0=1, \ell_1=0, \ell_2=1\end{cases}$
    }
  & Thm~\ref{thm18} &
     \\\hline
     $DD$ &  New & \makecell[l]{
       1, 1, 2, 1, 4, 5, 11, 11, 27, ...\\
       \text{General formula remains to be found}
       }
 &  Thm~\ref{thm21}  & Algebraic \\ \hline
\end{tabular}}
  \end{center}
\caption{Summary of results, generating functions are given in the corresponding theorems.}
\label{tt}
\end{table}

\section{Patterns of length one}
Our goal in this
section is to provide generating functions that count $p$-equivalence
classes with respect to the path length, whenever the pattern $p$ is
$U,$ $D,$ or $C$.  For each pattern, the method consists in
constructing a set of representative elements of $p$-equivalence
classes and enumerating them
with respect to the length by giving a generating function. From this,
we deduce a closed-form expression (that depends on the length) or a
recurrence relation for the number of corresponding $p$-equivalence
classes.

\subsection{Pattern $U$}
Let us define the set $\mathcal{A}$ consisting of the union of the set
$\mathcal{D}$ of Dyck paths with the set $\mathcal{D}'$ of paths in
$\mathcal{E}$ having only one catastrophe $C_k$, $k\geq 2$, located at
the end of the path. We will prove that there is a one-to-one
correspondence between $\mathcal{A}$ and the $U$-equivalence classes
of $\mathcal{E}$, and we enumerate them by providing a generating function that
counts $\mathcal{A}$ with respect to the length.

    \begin{thm}
    There is a bijection between $\mathcal{A}$ and the set of $U$-equivalence classes of $\mathcal{E}$.
    \label{thm1}
    \end{thm}

\noindent {\it Proof.} First, we prove that for every $P\in\mathcal{E}$ there is $Q\in \mathcal{A}$ such that $P$ and $Q$ belong to the same $U$-equivalence class. If we have $P\in\mathcal{D}$, then we choose $Q=P\in\mathcal{A}$. Otherwise, we  decompose $P=\alpha_1C_{k_1}\alpha_2C_{k_2}\ldots \alpha_rC_{k_r} \alpha_{r+1}$ where $r\geq 1$, $k_i\geq 2$ for $1\leq i\leq r$ and all $\alpha_i$ do not contain any catastrophe steps.  If $\alpha_{r+1}$ is empty, then we set $Q=\alpha_1D\alpha_2D\ldots \alpha_rC_{k}$ where $k=(k_1-1)+(k_2-1)+\ldots+ (k_{r-1}-1)+k_r=\sum_{i=1}^rk_i - (r-1)\geq 2r-(r-1)\geq 2$.
If $\alpha_{r+1}$ is not empty, it can be written $\alpha'D$, and we set $Q=\alpha_1D\alpha_2D\ldots \alpha_rD\alpha'C_k$ where $k=k_1-1+k_2-1+ k_{r-1}-1+k_r-1+1=\sum_{i=1}^rk_i-r+1\geq 2r-r+1\geq 2$. In these last two cases, we have $Q\in\mathcal{A}$ so that $P$ and $Q$  belong to the same $U$-equivalence class.

Now, let us prove that if $P$ and $Q$ are two paths of the same length
in $\mathcal{A}$ lying in the same $U$-equivalence class, then
$P=Q$. According to the decomposition of a Dyck path with catastrophes
(see the Introduction), we write either $P=\alpha$ or $P=\alpha
U\alpha_1 U\alpha_2 \ldots U\alpha_kC_k$ (respectively $Q=\alpha'$ or
$Q=\alpha'U\alpha'_1 U\alpha'_2 \ldots U\alpha'_\ell C_\ell$) for some
$k,\ell\geq 2$ and where $\alpha$, $\alpha'$, $\alpha_i$, $\alpha'_i$
are some Dyck paths. Since a Dyck path is characterized by the
positions of its up steps, and $P$ and $Q$ are in the same
$U$-equivalence class, we necessarily have $\alpha=\alpha'$,
$\alpha'_i=\alpha'_i$ for $i\leq \min\{k,\ell\}$, which implies that
$P=Q$. \hfill $\Box$

\begin{thm} The o.g.f.\ with respect to the length for the set $\mathcal{A}$ is given by
\begin{equation*}
    A(x)=\frac{\left(1-\sqrt{1-4 x^{2}}\right) \left(1-x\right)}{x \left(2x-1+\sqrt{1-4 x^{2}} \right)}.
\end{equation*}
 The series expansion of $A(x)$ is
      $$1+x^2+x^3+3x^4+4x^5+10x^6+15x^7+35x^8+56x^9+126x^{10}+210x^{11} +O\left(x^{12}\right).$$
      \label{thm2}
      \end{thm}
\noindent {\it Proof.}
A path in $\mathcal{A}$ is either a Dyck path in $\mathcal{D}$,  or a path  $\alpha U\alpha_1U\alpha_2U\alpha_3 \ldots U\alpha_kC_k$ where $k\geq 2$ and $\alpha,\alpha_1,\ldots, \alpha_k$ are some Dyck paths.  We deduce that
$A(x)=D(x)+D(x)\cdot\frac{x^2D(x)^2}{1-xD(x)}\cdot x$ where $D(x)=\frac{1-\sqrt{1-4 x^{2}}}{2 x^{2}}$ is the o.g.f.\ that counts Dyck paths with respect to the length (see~\cite{Deut}).
\hfill $\Box$

\begin{rem} With Theorem~\ref{thm2}, it is easy to check that  $x^2-2x+1+(-2x^2+3x-1)A(x)+(2x^2-x)A(x)^2=~0,$
and
the coefficient $a_n$ of $x^n$ in $A(x)$ is  ${n-1 \choose {\lfloor\frac{n-2}{2}\rfloor}}$ which generates a shift of the sequence \oeis{A037952} in OEIS~\cite{OEIS}. Using a classical analysis~\cite{flajolet,orlov} of dominant singularity of $A(x)$, $a_n$ has the asymptotic approximation $\frac{2^n}{\sqrt{2\pi n}}.$
\end{rem}

\subsection{Pattern $D$}
Let us define the set $\mathcal{B}$ consisting of the union of the set $\mathcal{D}$ of Dyck paths with the set $\mathcal{D}''$  of paths $P\in \mathcal{E}$ having only one catastrophe.

    \begin{thm}
    There is a bijection between $\mathcal{B}$ and the set of $D$-equivalence classes of $\mathcal{E}$.
    \label{thm3}
    \end{thm}

\noindent {\it Proof.} First, we prove that for every $P\in\mathcal{E}$ there is $Q\in \mathcal{B}$ such that $P$ and $Q$ belong to the same $D$-equivalence class. If we have $P\in\mathcal{D}$, then we choose $Q=P\in\mathcal{B}$;
 otherwise, we  decompose $P=\alpha_1C_{k_1}\alpha_2C_{k_2}\ldots \alpha_rC_{k_r} \alpha_{r+1}$ where $r\geq 1$, $k_i\geq 2$ for $1\leq i\leq r$ and all $\alpha_i$ do not contain any catastrophe steps.  If $\alpha_{r+1}$ is empty, then we set $Q=\alpha_1U\alpha_2U\ldots \alpha_rC_{k}$ where  $k=(k_1+1)+(k_2+1)+\ldots + (k_{r-1}+1)+k_r=r-1+\sum_{i=1}^rk_i$.  If $\alpha_{r+1}$ is not empty, it is necessarily a Dyck path, and we set $Q=\alpha_1U\alpha_2U\ldots \alpha_rC_k\alpha_{r+1}$ where $k=k_1+1+k_2+1+ k_{r-1}+1+k_r=r-1+\sum_{i=1}^rk_i$. In all these cases, we have $Q\in\mathcal{B}$ so that $P$ and $Q$  belong to the same $D$-equivalence class.

With a similar argument as for the proof of Theorem~\ref{thm1}, it is easy to prove
that if $P$ and $Q$ are two paths in $\mathcal{B}$ lying in the same class, then $P=Q$.

\hfill $\Box$
\begin{thm} The o.g.f.\ with respect to the length for the set $\mathcal{B}$ is given by
\begin{equation*}
    B(x)=\frac{\left(1-\sqrt{1-4 x^{2}}\right) \left(2 x^{3}-4 x^{2}+1+(2 x^{2}-1) \sqrt{1-4 x^{2}}\right)}{2 x^{4} \left(2x-1+\sqrt{1-4 x^{2}} \right)}
.
\end{equation*}
 The series expansion of $B(x)$ is
      $$1+x^2+x^3+3x^4+5x^5+11x^6+21x^7+42x^8+84x^9+162x^{10}+ 330x^{11} +O\left(x^{12}\right).$$
      \label{thm4}
      \end{thm}
\noindent {\it Proof.}
A path in $\mathcal{B}$ is either a Dyck path in $\mathcal{D}$,  or a path  of the form $\alpha U\alpha_1U\alpha_2U\alpha_3 \ldots U\alpha_kC_k\beta$ where $k\geq 2$ and $\alpha,\alpha_1,\ldots, \alpha_k,\beta$ are some Dyck paths.  We deduce that $B(x)$ satisfies the functional equation $B(x)=D(x)+\frac{x^3D(x)^4}{1-xD(x)}$ where $D(x)=\frac{1-\sqrt{1-4 x^{2}}}{2 x^{2}}$ is the o.g.f.\ that counts Dyck paths with respect to the length (see~\cite{Deut}).
\hfill $\Box$

\begin{rem} We have $5x^3-4x^2-x+1+(-2x^4-5x^3+5x^2+x-1)B(x)+(2x^5-x^4)B(x)^2=0,$ and  the coefficient $b_n$ of $x^n$ in $B(x)$ is
  $$b_n=\begin{cases}
  \displaystyle
      {n \choose \frac{n-3}{2}} &\text{ if } n \text{ is odd,} \vspace{1ex} \\
        \displaystyle
{n \choose  \frac{n-4}{2} }+
{n \choose \frac{n}{2}}\slash\left(\frac{n}{2}+1\right) &\text{ otherwise.}
\end{cases}$$
This last result can be easily obtained by combining already known
formul\ae\ (see \oeis{A002054} and \oeis{A344191} in~\cite{OEIS}), or directly using
Maple. The sequence $b_n$ has an asymptotic approximation
$\frac{2^{n+\frac{1}{2}}}{\sqrt{\pi n}}.$
\end{rem}

\subsection{Pattern $C$}

In this section, two paths of the same length are $C${\it -equivalent}
whenever, for any $k\geq 2$ they have the same positions of the
occurrences of $C_k$. For instance,
$$UUDDUUC_2UUUC_3UUC_2UD$$
is $C$-equivalent to
$$UDUDUUC_2UUUC_3UUC_2UD,$$
but not to
$$UUUDUUC_4UUUC_3UUC_2UD.$$

Let $\mathcal{C}$ be the set of paths $P\in\mathcal{E}$ such that:

\noindent ($i$)  $P=(UD)^k$ for $k\geq 0$, or

\noindent ($ii$) $P=(UD)^{\ell_1}U^{k_1}C_{k_1}  (UD)^{\ell_2}U^{k_2}C_{k_2} \ldots (UD)^{\ell_r}U^{k_r}C_{k_r} (UD)^{\ell_{r+1}}$ with $r\geq~1$, $\ell_i\geq  0$ for $1\leq i\leq r+1$, and $k_i\geq 2$ for $1\leq i\leq r$.

   \begin{thm}
    There is a bijection between $\mathcal{C}$ and the set of $C$-equivalence classes of $\mathcal{E}$.
    \label{thm5}
    \end{thm}

\noindent {\it Proof.} First, we prove that for every $P\in\mathcal{E}$ there is $Q\in \mathcal{C}$ such that $P$ and $Q$ belong to the same $C$-equivalence class. If we have $P\in\mathcal{D}$, then we choose $Q=(UD)^k\in\mathcal{C}$ with $k=|P|\slash 2$. Otherwise, we  decompose $P=\alpha_1C_{k_1}\alpha_2C_{k_2}\ldots \alpha_rC_{k_r} \alpha_{r+1}$ where $r\geq 1$, $k_i\geq 2$ for $1\leq i\leq r$ and all $\alpha_i$ do not contain any catastrophe steps.  We set $Q=(UD)^{j_1}U^{k_1}C_{k_1} (UD)^{j_2}U^{k_2}C_{k_2}\ldots \break (UD)^{j_r}U^{k_r}C_{k_r}(UD)^{j_{r+1}}$ where $j_i+k_i=|\alpha_i|$, $1\leq i\leq r$ and $j_{r+1}=|\alpha_{r+1}|$.
We have $Q\in\mathcal{C}$ so that $P$ and $Q$  belong to the same $C$-equivalence class. Due to the form of $Q$, it is straightforward to see that if $P$ and $Q$ are two paths in $\mathcal{C}$ lying in the same  class, then $P=Q$. \hfill $\Box$

\begin{thm} The o.g.f.\ with respect to the length for the set $\mathcal{C}$ is given by
\begin{equation*}
    C(x)=\frac{1-x}{1-x-x^2}.
\end{equation*}
 The series expansion of $C(x)$ is
      $$1+x^2+x^3+2x^4+3x^5+5x^6+8x^7+13x^8+21x^9+34x^{10}+55x^{11} +O\left(x^{12}\right).$$
      \label{thm6}
      \end{thm}
\noindent {\it Proof.}
Due to the definition of the set $\mathcal{C}$, we obtain directly the functional equation
$C(x)=\frac{1}{1-x^2} + \frac{1}{1-x^2}\cdot\frac{C'(x)}{1-C'(x)}$ where $C'(x)=\frac{x^3}{1-x}\cdot\frac{1}{1-x^2}$  is the o.g.f.\ that counts the paths of the form
 $(UD)^{\ell}U^{k}C_{k}$ for $\ell\geq 0$ and $k\geq 2$.
\hfill $\Box$

\begin{rem} From Theorem~\ref{thm6}, we deduce that the coefficient $c_n$ of $x^n$ in
 $C(x)$ satisfies $c_{n}=c_{n-1}+c_{n-2}$ for $n\geq 2$ with $c_0=1$ and $c_1=0$, which generates the sequence \oeis{A212804} in~\cite{OEIS} (it is a variant of the  well known Fibonacci sequence A000045). Using the classical method for asymptotic approximation (see ~\cite{flajolet,orlov}), we have
 $$c_n\sim\frac{3-\sqrt{5}}{5-\sqrt{5}}\cdot \left( \frac{\sqrt{5}-1}{2}\right)^n.$$
\end{rem}


\section{Patterns of length two}

In this section, we provide generating functions that count
$p$-equivalence classes with respect to the length whenever the
pattern $p$ is $UU$, $UD$, $UC$, $DC$, $CU$, $DU$ or $DD$. The 
techniques used are similar, although more elaborate, to those used in
  the previous section.

\subsection{Pattern $UU$}

Let $\mathcal{F}$ be the set of paths $P\in\mathcal{E}$ such that:

\noindent ($i$) $P=(UD)^k$ for $k\geq 0$, or

\noindent ($ii$) $P=(UD)^{\ell_0}U^{k_1} \alpha_1 U^{k_2}\alpha_2 \ldots U^{k_r}\alpha_r$ where $r\geq 1$, $\ell_0\geq 0$,  $ k_i\geq 2$ for $1\leq i\leq r$, all $\alpha_{i}$ for  $1\leq i\leq r-1$  are either $(DU)^kD$ or $(DU)^kDD$ for some $k\geq 0$, and $\alpha_{r}$ is either $(DU)^kC_s$ or $(DU)^kDC_s$ for some $k\geq 0$ and with $s\geq 1$ is so that the path ends on the $x$-axis (note that $s$ can be 1, and in this case $C_1=D$).

 \begin{thm}
    There is a bijection between $\mathcal{F}$ and the set of $UU$-equivalence classes of $\mathcal{E}$.
    \label{thm7}
    \end{thm}

\noindent {\it Proof.} First, we prove that for every $P\in\mathcal{E}$ there is $Q\in \mathcal{F}$ such that $P$ and $Q$ belong to the same $UU$-equivalence class. If $P$ does not contain occurrences of $UU$, then we choose $Q=(UD)^k\in\mathcal{F}$ where $k=|P|\slash 2$. Otherwise, we decompose $P=\alpha_0U^{k_1}\alpha_1U^{k_2} \ldots \alpha_{r-1}U^{k_r}\alpha_{r}$ where $r\geq 1$, $k_i\geq 2$ for $1\leq i\leq r$, and such that all occurrences of $UU$ in $P$ belong necessarily to a run $U^{k_i}$ for some $i$.

We set $Q=(UD)^{\ell_0}U^{k_1} \beta_1 U^{k_2} \ldots \beta_{r-1}U^{k_r}\beta_{r}$ where $\ell_0=|\alpha_0|$, and for $1\leq i\leq r-1$,
$\beta_i=(DU)^{\frac{t_i-1}{2}}D$ if $t_i=|\alpha_i|$ is odd,
$\beta_i=(DU)^{\frac{t_i-2}{2}}DD$ otherwise; finally, we set  $\beta_{r}=(DU)^{\frac{t_{r}-1}{2}}C_s$ if $t_{r}=|\alpha_{r}|$ is odd,
and otherwise,  $\beta_{r}=(DU)^{\frac{t_{r}-2}{2}}DC_s$ with $s\geq 1$ is chosen so that $Q$ ends on the $x$-axis.

We have $Q\in\mathcal{F}$ so that $P$ and $Q$  belong to the same $UU$-equivalence class. It is straightforward to see that if $P$ and $Q$ are two paths in $\mathcal{F}$ lying in the same class, then $P=Q$.
 \hfill $\Box$

\begin{thm} The o.g.f.\ with respect to the length for the set $\mathcal{F}$ is given by
\begin{equation*}
    F(x)={\frac {x-1}{ \left( x+1 \right)  \left( {x}^{3}-{x}^{2}+2\,x-1
 \right) }}.
\end{equation*}
 The series expansion of $F(x)$ is
      $$1+x^2+x^3+3x^4+4x^5+8x^6+13x^7+24x^8+41x^9+73x^{10}+127x^{11} +O\left(x^{12}\right).$$
      \label{thm8}
      \end{thm}
\noindent {\it Proof.} Due to the definition of $\mathcal{F}$, we have $F(x)=\frac{1}{1-x^2}+\frac{1}{1-x^2}\cdot \frac{1}{1-\frac{x^3}{(1-x)^2}}\frac{x^3}{(1-x)^2}$.
A simple calculation completes the proof.
  \hfill $\Box$
 \begin{rem} The coefficient $f_n$ of $x^n$ in
 $F(x)$ satisfies $f_{n}=f_{n-1}+f_{n-2}+f_{n-4}$ with $f_0=1, f_1=0, f_2=1, f_3=1$, which generates the sequence \oeis{A347493} in~\cite{OEIS}. An asymptotic approximation of $f_n$ is $$f_n\sim \frac{(a-1)(a^2-a+2)^n}{a(a+1)(-3a^2+2a-2)}\approx 0.26212\cdot 1.75487^n,$$ where $a=\frac{(44+12\sqrt{69})^{2/3}+2(44+12\sqrt{69})^{1/3}-20}{6(44+12\sqrt{69})^{1/3}}$.
\end{rem}

 \subsection{Pattern $UD$}

Let $\mathcal{G}$ be the set of paths $P\in\mathcal{E}$ such that either :

\noindent ($i$) $P=(UD)^k$ for some $k\geq 0$, or

\noindent ($ii$) $P=(UD)^\ell U^kC_k (UD)^m$ with $\ell,m\geq 0$ and $k\geq 2$, or

\noindent ($iii$) $P=(UD)^{\ell_{0}}U^{k_1}(UD)^{\ell_1}U^{k_2}(UD)^{\ell_2}\ldots U^{k_r}(UD)^{\ell_r}U^{k_{r+1}}C_s(UD)^{\ell_{r+1}}$ with $r\geq 1$, $\ell_{0}, \ell_{r+1}\geq 0$, $k_i\geq 1$ for $1\leq i\leq r$, $k_{r+1}\geq 0$,  $\ell_i\geq 1$ for $1\leq i\leq r$, and $s\geq 1$ is so that the path ends on the $x$-axis (note that $s$ can be one).

 \begin{thm}
    There is a bijection between $\mathcal{G}$ and the set of $UD$-equivalence classes of $\mathcal{E}$.
    \label{thm9}
    \end{thm}

\noindent {\it Proof.} First, we prove that for every $P\in\mathcal{E}$ there is $Q\in \mathcal{G}$ such that $P$ and $Q$ belong to the same $UD$-equivalence class. If $P$ satisfies the case ($i$), then we obviously set $Q=P$. If $P=(UD)^\ell \alpha (UD)^m$ with $l,m\geq 0$ and $\alpha$ is a nonempty path in $\mathcal{E}$ avoiding $UD$, then we set $Q=(UD)^\ell U^kC_k (UD)^m$. Otherwise, we decompose $P=(UD)^{\ell_{0}}\alpha_1(UD)^{\ell_1}\alpha_2(UD)^{\ell_2}\ldots \alpha_r(UD)^{\ell_r}\alpha_{r+1}$ with all $\alpha_i$ being non-empty partial paths avoiding $UD$ (except $\alpha_{r+1}$ that can be empty). If $|\alpha_{r+1}|>1$ then we set
$Q=(UD)^{\ell_{0}}U^{k_1}(UD)^{\ell_1}U^{k_2}(UD)^{\ell_2}\ldots \\ U^{k_r}(UD)^{\ell_r}U^{k_{r+1}}C_s$ where $k_i=|\alpha_i|$ and $k_{r+1}+1=|\alpha_{r+1}|$; otherwise, if $|\alpha_{r+1}|=0$ then we set
$Q=(UD)^{\ell_{0}}U^{k_1}(UD)^{\ell_1}U^{k_2}(UD)^{\ell_2}\ldots U^{k_r-1}C_s(UD)^{\ell_r}$  where $k_i=|\alpha_i|$ and $s$ is so that $Q$ ends on the $x$-axis.
For all these cases, we have $Q\in\mathcal{G}$ so that $P$ and $Q$ belong to the same $UD$-equivalence class. Due to the form of a path in $\mathcal{G}$, if $P$ and $Q$ are two paths in $\mathcal{G}$ lying in the same $UD$-equivalence class, the $P=Q$.\hfill $\Box$

\begin{thm} The o.g.f.\ with respect to the length for the set $\mathcal{G}$ is given by
\begin{equation*}
    G(x)=\frac{2 x^{3}-2 x +1}{\left(1-x^{2}-x \right) \left(1-x \right)^{2} \left(x +1\right)}.
\end{equation*}
 The series expansion of $G(x)$ is
      $$1+x^2+x^3+3x^4+5x^5+10x^6+17x^7+30x^{8}+50x^{9}+84x^{10}+138x^{11}
 +O\left(x^{12}\right).$$
 \label{thm10}
      \end{thm}
\noindent {\it Proof.} According to the different cases in the definition of $\mathcal{G}$,
$G(x)$ is given by: $$1+\frac{x^2}{1-x^2}+\frac{1}{1-x^2}\frac{x^3}{1-x}\frac{1}{1-x^2} +\frac{1}{1-x^2} \frac{R(x)}{1-R(x)}\frac{x}{1-x}\frac{1}{1-x^2}$$ where $R(x)=\frac{x^3}{(1-x)(1-x^2)}$ is the o.g.f.\ that counts the paths of the form $U^{k}(UD)^{\ell}$ with $k,\ell\geq 1$.
  \hfill $\Box$
 \begin{rem} The coefficient $g_n$ of $x^n$ in
 $G(x)$ satisfies $g_n = 2g_{n-1} + g_{n-2} - 3g_{n-3} + g_{n-5}$ for $n>4$ with $g_1=0$, $g_0=g_2=g_3=1$, $g_4=3$, which generates the sequence \oeis{A215004} in~\cite{OEIS}. An asymptotic approximation  of $g_n$
  is $$g_n\sim\frac{4(5-2\sqrt{5})}{5(\sqrt{5}-3)^2}\cdot \left(\frac{\sqrt{5}+1}{2}\right)^n\approx 0.72360\cdot 1.61803^n.$$

\end{rem}

 \subsection{Pattern $UC$}

In this section, two paths are $UC$-equivalent whenever they coincide on all their occurrences  $UC_k$ for  $k\geq 2$. For instance, $UUDUC_2$ and $UDUUC_2$ are $UC$-equivalent, while $UDUUC_2$ and $UUUUC_4$ are not.

Let $\mathcal{I}_1$ be the set of paths of length $n\geq 0$, $n\notin\{1,3\}$, defined by either

\noindent ($i$) $(UD)^{\frac{n}{2}}$ if $n$ is even, or

\noindent ($ii$) $(UD)^{\frac{n-5}{2}}UUUDC_2$ if is $n$ is odd.

\noindent Let $\mathcal{I}$ be the set consisting of the union of $\mathcal{I}_1$ with the set of paths of the form
$\alpha_1 U^{k_1}C_{k_1} \alpha_2 U^{k_2}C_{k_2} \ldots \alpha_r U^{k_r}C_{k_r}\alpha_{r+1}$
where $r\geq 1$, all values $k_i$ are at least $2$, and all $\alpha_i$ are in $\mathcal{I}_1$.

 \begin{thm}
    There is a bijection between $\mathcal{I}$ and the set of $UC$-equivalence classes of $\mathcal{E}$.
    \label{thm11}
    \end{thm}

\noindent {\it Proof.}  First, we prove that for every $P\in\mathcal{E}$ there is $Q\in \mathcal{I}$ such that $P$ and $Q$ belong to the same $UC$-equivalence class. If $P$  does not contain any pattern $UC$, then we set $Q=(UD)^{\frac{n}{2}}$ if $n$ is even, and $Q=(UD)^{\frac{n-5}{2}}UUUDC_2$ otherwise. If $P$ contains $r\geq 1$ occurrences of $UC$, then we decompose $P=\alpha_1UC\alpha_2UC\ldots \alpha_rUC\alpha_{r+1}$ where all $\alpha_i$ are paths avoiding $UC$. Showing the size of every catastrophes, we obtain the decomposition
$P=\alpha_1UC_{k_1}\alpha_2UC_{k_2}\ldots \alpha_rUC_{k_r}\alpha_{r+1}$ where $k_i\geq 2$ for $1\leq i\leq r$.

We set  $Q=\beta_1U^{k_1}C_{k_1}\beta_2U^{k_2}C_{k_2}\ldots \beta_rU^{k_r}C_{k_r}\beta_{r+1}$ where all $\beta_i$ are in $\mathcal{I}_1$ (all $\beta_i$ are entirely determined by the length of each $\alpha_i$).
We have $Q\in\mathcal{I}$ so that $P$ and $Q$ belong to the same $UC$-equivalence class. Obviously, due to the definition of $\mathcal{I}$, if $P$ and $Q$ belong to $\mathcal{I}$ in the same $UC$-equivalence class, then $P=Q$.\hfill $\Box$

\begin{thm} The o.g.f.\ with respect to the length for the set $\mathcal{I}$ is given by
\begin{equation*}
    I(x)={\frac {{x}^{5}-2\,{x}^{4}+2\,{x}^{3}-2\,{x}^{2}+2\,x-1}{{x}^{7}-{x}^{
6}+{x}^{5}-{x}^{4}+{x}^{3}-{x}^{2}+2\,x-1}}.
\end{equation*}
 The series expansion of $I(x)$ is
      $$1+x^2+x^3+2x^4+4x^5+5x^6+9x^7+15x^{8}+24x^{9}+40x^{10}+65x^{11}
 +O\left(x^{12}\right).$$
 \label{thm12}
      \end{thm}
\noindent {\it Proof.}
The o.g.f.\ for $\mathcal{I}_1$ is $I_1(x)=\frac{1}{1-x}-x-x^3$  and the o.g.f.\ for $\mathcal{I}$ is $I_1(x)+I_1(x)\cdot \frac{I_2(x)}{1-I_2(x)}$ where $I_2(x)=\frac{I_1(x)x^3}{1-x}$ is the o.g.f.\ of a nonempty sequence of terms of the form  $\alpha U^{k}C_{k}$ with
$\alpha\in \mathcal{I}_1$ and $k\geq 2$.
\hfill $\Box$

 \begin{rem} The coefficient $i_n$ of $x^n$ in
 $I(x)$ satisfies $i_n =2i_{n-1}-i_{n-2}+i_{n-3}-i_{n-4}+i_{n-5}-i_{n-6}+i_{n-7}$ for $n\geq 7$ with $i_0=1, i_1=0, i_2=1, i_3=1, i_4=2, i_5=4$, and $i_6=5$. This sequence does not appear in~\cite{OEIS}. An asymptotic approximation of $i_n$ is $$ i_n\sim  0.28134\cdot 1.64072^n.$$
\end{rem}

 \subsection{Pattern $DC$}
In this section, two paths are $DC$-equivalent whenever they coincide on all their occurrences  $DC_k$ for  $k\geq 2$.

Let $\mathcal{J}_1$ be the set of paths of length $n\geq 0$, $n\neq 1$, defined by either

\noindent ($i$) $(UD)^{\frac{n}{2}}$ if $n$ is even, or

\noindent ($ii$) $(UD)^{\frac{n-3}{2}}UUC_2$ if is $n$ is odd.

\noindent Let $\mathcal{J}$ be the set consisting of the union of $\mathcal{J}_1$ with the set of paths of the form
$\alpha_1 U^{k_1+1}DC_{k_1} \alpha_2 U^{k_2+1}DC_{k_2} \ldots \alpha_r U^{k_r+1}DC_{k_r}\alpha_{r+1}$
where $r\geq 1$, all values $k_i$ are at least $2$, and all $\alpha_i$ are in $\mathcal{J}_1$.

 \begin{thm}
    There is a bijection between $\mathcal{J}$ and the set of $DC$-equivalence classes of $\mathcal{E}$.
    \label{thm13}
    \end{thm}

\noindent {\it Proof.} The  proof is obtained {\it mutatis mutandis} as for Theorem~\ref{thm11}.
\hfill $\Box$

\begin{thm} The o.g.f.\ with respect to the length for the set $\mathcal{J}$ is given by
\begin{equation*}
    J(x)=\frac{x^{3}-2 x^{2}+2 x -1}{x^{7}-x^{6}+x^{5}-x^{2}+2 x -1}.
\end{equation*}
 The series expansion of $J(x)$ is
      $$1+x^2+x^3+x^4+2x^5+2x^6+4x^7+6x^{8}+9x^{9}+14x^{10}+20x^{11}+O\left(x^{12}\right).$$
      \label{thm14}
      \end{thm}
\noindent {\it Proof.} The o.g.f.\ for $\mathcal{J}_1$ is $J_1(x)=\frac{1}{1-x}-x$  and the o.g.f.\ for $\mathcal{J}$ is $J_1(x)+J_1(x)\cdot \frac{J_2(x)}{1-J_2(x)}$ where $J_2(x)=\frac{J_1(x)x^5}{1-x}$ is the o.g.f.\ of a nonempty sequence of terms of the form  $\alpha U^{k+1}DC_{k}$ with
$\alpha\in \mathcal{J}_1$ and $k\geq 2$. \hfill $\Box$

 \begin{rem} The coefficient $j_n$ of $x^n$ in
 $J(x)$ satisfies $j_n=2j_{n-1}-j_{n-2}+j_{n-5}-j_{n-6}+j_{n-7}$ for $n\geq 7$ with $j_0=1, j_1=0, j_2=j_3=j_4=1, j_5=2, j_6=2$ which is not in~\cite{OEIS}. An asymptotic approximation of $j_n$ is $$j_n\sim  0.25317\cdot 1.48698^n.$$
\end{rem}

\subsection{Pattern $CU$}

In this section, two paths are $CU$-equivalent whenever they coincide on all their occurrences  $C_kU$  for  $k\geq 2$.

Let $\mathcal{K}_1$ be the set of paths of even length defined by $(UD)^k$  for $k\geq 0$.

\noindent Let $\mathcal{K}_2$ be the set of paths of length $n\neq 1$ defined by $(UD)^k$, $k \geq 1$, if $n$ is even,  and $(UD)^kUUC_2$, $k\geq 0$,  if $n$ is odd.

\noindent Let $\mathcal{K}$ be the set consisting of the union of $\mathcal{K}_2$ with the set of paths of the form
$\beta_1 U^{k_1} C_{k_1}
\beta_2 U^{k_2} C_{k_2}  \ldots
\beta_r U^{k_r} C_{k_r}
\beta_{r+1}$ where $k_i\geq 2$ for $1 \le i < r $, and for $1 \le i<r+1$, $\beta_i \in \mathcal{K}_1$,
    and $\beta_{r+1}$ is a path from $\mathcal{K}_2$.

 \begin{thm}
    There is a bijection between $\mathcal{K}$ and the set of $CU$-equivalence classes of $\mathcal{E}$.
    \label{thm15}
    \end{thm}

\noindent {\it Proof.} First, we prove that for every $P\in\mathcal{E}$ there is $Q\in \mathcal{K}$ such that $P$ and $Q$ belong to the same $CU$-equivalence class. If $P$  does not contain any pattern $CU$, then we set $Q=(UD)^{\frac{n}{2}}$ whenever $n=|P|$ is even, and we set $Q=(UD)^{\frac{n-3}{2}}UUC_2$ otherwise; in both cases $Q\in\mathcal{K}_2$. If $P$ contains $r\geq 1$ occurrences of $CU$, then we decompose
$P=\alpha_1CU\alpha_2CU\ldots \alpha_rCU\alpha_{r+1}$ where all $\alpha_i$ are paths avoiding $CU$. Showing the size of every catastrophe,
we obtain  the decomposition
$P=\alpha_1C_{k_1}U\alpha_2C_{k_2}U\ldots \alpha_rC_{k_r}U\alpha_{r+1}$ where all $k_i$ are at least $2$.
So, we set
$Q=
\beta_1 U^{k_1} C_{k_1}
\beta_2 U^{k_2} C_{k_2}  \ldots
\beta_r U^{k_r} C_{k_r}
\beta_{r+1}$ where
$|\alpha_1| = |\beta_1|$,
$|\alpha_i|+1 = |\beta_i| + k_i$
for $1 < i < r + 1$,
$|\alpha_{r+1}|+1 = |\beta_{r+1}|$, and

\begin{itemize}
    \item for $1 \le i<r+1$, $\beta_i \in \mathcal{K}_1$;
    \item $\beta_{r+1}$ is a path from $\mathcal{K}_2$.
\end{itemize}

We have $Q\in\mathcal{K}$ so that $P$ and $Q$ belong to the same $CU$-equivalence class. Obviously, due to the definition of $\mathcal{K}$, if $P$ and $Q$ belong to $\mathcal{K}$ in the same $CU$-equivalence class, then $P=Q$.\hfill $\Box$

\begin{thm} The o.g.f.\ with respect to the length for the set $\mathcal{K}$ is given by
\begin{equation*}
    K(x)=\frac{x^{4}+x -1}{x^{2}+x -1}.
\end{equation*}
 The series expansion of $K(x)$ is
      $$1+x^2+x^3+x^4+2x^5+3x^6+5x^7+8x^8+13x^9+21x^{10}+34x^{11}
 +O\left(x^{12}\right).$$
 \label{thm16}
      \end{thm}
\noindent {\it Proof.}
 Due to the definition of $\mathcal{K}$, the o.g.f.\ $K(x)$ is given by
 $$K(x)=\frac{1}{1-x}-x+ \frac{1}{1-x^2} \cdot\frac{1}{1-\frac{x^3}{(1-x)(1-x^2)}}
 \cdot \frac{x^3}{1-x}\cdot\frac{x^2}{1-x}.$$

\hfill $\Box$

 \begin{rem} The coefficient $k_n$ of $x^n$ in
 $K(x)$ satisfies $k_n=k_{n-1}+k_{n-2}$ for $n\geq 5$ with $k_0=1, k_1=0, k_2=k_3=k_4=1$ which is a shift of the Fibonacci sequence \oeis{A000045} in~\cite{OEIS}. An asymptotic approximation of $k_n$ is
 $$k_n\sim \frac{2(\sqrt{5}-2)}{5-\sqrt{5}}\cdot \left(\frac{2}{\sqrt{5}-1}\right)^n\approx 0.17082\cdot 1.61803^n.$$
\end{rem}

\subsection{Pattern $DU$}
Let $\mathcal{L}$ be set of paths $P\in\mathcal{E}$ such that $$P=U^{\ell_1}(DU)^{k_1}U^{\ell_2}(DU)^{k_2}\ldots U^{k_r} (DU)^{k_r}U^{\ell_{r+1}}C_s$$ with $r\geq 1$, $\ell_i,k_i\geq 1$ for $1\leq i\leq r$, $\ell_{r+1}\geq 0$, and  $s$ is so that the path ends on the $x$-axis (note that $s$ can be $1$, and in this case $C_1=D$).

 \begin{thm}
    There is a bijection between $\mathcal{L}$ and the set of $DU$-equivalence classes of $\mathcal{E}$.
    \label{thm17}
    \end{thm}

\noindent {\it Proof.} First, we prove that for every $P\in\mathcal{E}$ there is $Q\in \mathcal{G}$ such that $P$ and $Q$ belong to the same $DU$-equivalence class. We decompose $P=\alpha_1(DU)^{k_1}\alpha_2(DU)^{k_2}\alpha_3(DU)^{k_3}\ldots \alpha_r(DU)^{k_r}\alpha_{r+1}$ with all $\alpha_i$ being non-empty paths avoiding $DU$, with $r\geq 0$, and all $k_i$ are at least one. Then we set $Q=U^{\ell_1}(DU)^{k_1}U^{\ell_2}(DU)^{k_2}U^{\ell_3}(DU)^{k_3}\ldots U^{\ell_r}(DU)^{k_r}U^{\ell_{r+1}-1}C_s$ where $\ell_i=|\alpha_i|$, and where $s$ is so that the path ends on the $x$-axis.
For all these cases, we have $Q\in\mathcal{L}$ so that $P$ and $Q$ belong to the same $DU$-equivalence class. Due to the form of a path in $\mathcal{L}$, if $P$ and $Q$ are two paths in $\mathcal{L}$ lying in the same $DU$-equivalence class, then $P=Q$.\hfill $\Box$

\begin{thm} The o.g.f.\ with respect to the length for the set $\mathcal{L}$ is given by
\begin{equation*}
    L(x)=\frac{1-x}{1-x-x^2}.
\end{equation*}
 The series expansion of $L(x)$ is
      $$1+x^2+x^3+2x^4+3x^5+5x^6+8x^7+13x^{8}+21x^{9}+34x^{10}+55x^{11}
 +O\left(x^{12}\right).$$
 \label{thm18}
      \end{thm}
\noindent {\it Proof.} According to the definition of $\mathcal{L}$ and handling separately the cases where $r=1$ and $r\geq 2$, the o.g.f.\ $L(x)$ is given by: $$1+\frac{x^2}{1-x^2}+\frac{R(x)}{1-R(x)}\cdot \frac{x}{1-x}$$ where $R(x)=\frac{x^3}{(1-x)(1-x^2)}$ is the o.g.f.
that counts paths of the form $U^{\ell}(DU)^{k}$ with $k,\ell\geq 1$.
  \hfill $\Box$
 \begin{rem} The coefficient $\ell_n$ of $x^n$ in
 $G(x)$ satisfies $\ell_n = \ell_{n-1} + \ell_{n-2}$ for $n\geq 3$ with $\ell_0=1, \ell_1=0, \ell_2=1$, which generates the sequence \oeis{A212804} in~\cite{OEIS} (a variant of the well known Fibonacci sequence A000045). An asymptotic approximation of $\ell_n$ is $$\ell_n\sim\frac{(\sqrt{5}-3)}{5-\sqrt{5}}\cdot \left(\frac{2}{\sqrt{5}-1}\right)^n\approx 0.27639\cdot 1.61803^n.$$
\end{rem}

 \subsection{Pattern $DD$}
 \label{sDD}

 In~\cite{Bari1}, the authors exhibited a set $\overline{\mathcal{A}}$
 (denoted $\mathcal{A}$ in~\cite{Bari1}) of representatives for the
 $DD$-equivalence classes of Dyck paths: $\overline{\mathcal{A}}$
 consists of Dyck paths $P$ satisfying

 \medskip
 \noindent {\bf Condition $(C)$:}
 {\it $P$ avoids $UUDU$ and the minimal ordinate reached  by every occurrence of $UDU$ is
   at most one (i.e.~the height of every occurrence of $UDU$ is at most $1$).
 }
 \medskip

 \noindent For a Dyck path $P$, we denote by $\psi(P)$ the
 unique path in $\overline{\mathcal{A}}$ $DD$-equivalent to $P$. In~\cite{Bari1},
 it is shown how $\psi(P)$ can be constructed from $P$ after several
 transformations preserving the positions of all occurrences of $DD$.
 The following proposition extends this result for the set $\mathcal{D}'$ of Dyck paths
 with only one catastrophe $C_k$ at the end.

 \begin{prop}
 For any path $P\in \mathcal{D}'$ ending with catastrophe $C_k$,
 $k\geq 2$, there is a unique path in $\mathcal{D}'$, denoted
 $\psi(P)$, which is $DD$-equivalent to $P$, ending with catastrophe
 $C_k$ (the same $k$ as for $P$), and satisfying Condition~$(C)$.
 \label{prop1}
 \end{prop}

 \noindent {\it Proof.} Let $Q$ be the Dyck path obtained from $P$ by
 replacing $C_k$ with $D^k$. Using~\cite{Bari1}, there is a unique
 path $\psi(Q)$ in $\overline{\mathcal{A}}$ in the class of $Q$. Replacing the
 last $D^k$ of $\psi(Q)$ with $C_k$, there is a unique Dyck path with
 only catastrophe $C_k$ at the end equivalent to $P$.  \hfill $\Box$

\begin{figure}[H]
  \begin{center}
    \includegraphics[height=6em]{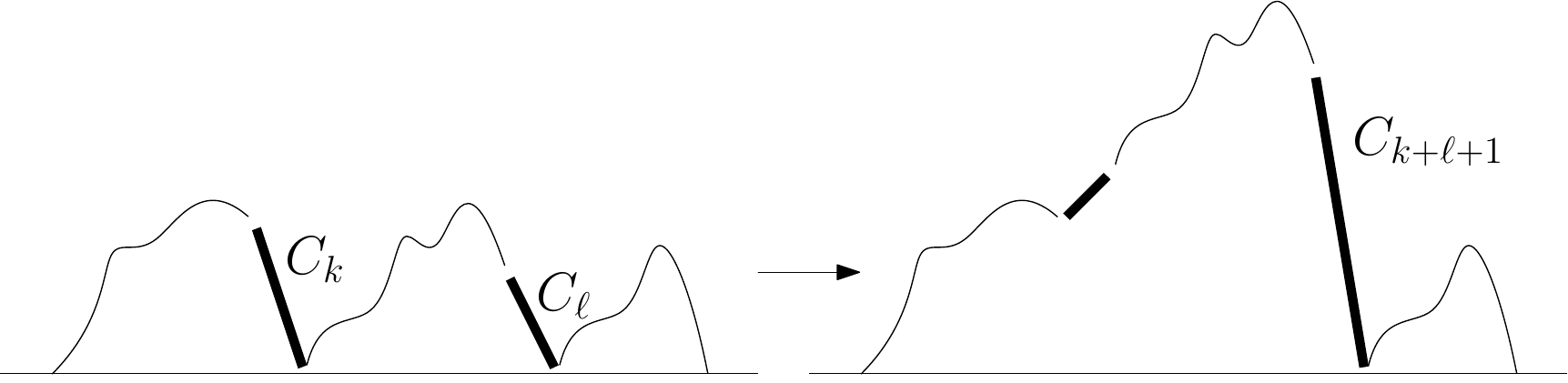}
  \end{center}
  \caption{Example for Lemma~\ref{dd_1} with two catastrophes.}
  \label{dd_1f}
\end{figure}

 \begin{lem}
   Let $P$ be a path in $\mathcal{E}$.  If $P$ contains at least two
   catastrophes, then there exists $Q\in \mathcal{E}$ in the class of
   $P$, such that $Q$ has only one catastrophe.
   \label{dd_1}
\end{lem}
\noindent {\it Proof.} We obtain $Q$ from $P$ by replacing with $U$
every catastrophe, except the last, and by increasing the size of the
last catastrophe so that $Q\in \mathcal{E}$. See Figure~\ref{dd_1f} for an example with two catastrophes.\hfill $\Box$

\medskip

Considering this previous lemma, we can focus our study on the set of
Dyck paths with at most one catastrophe.
 Recall that $\mathcal{D}''$ is the set of paths in $\mathcal{E}$ having only one catastrophe (see Section 2.2).

\begin{lem}
Let $P$ be a path of length $n$ in $\mathcal{D}''$ (i.e.~$P$ contains
only one catastrophe). The size of the catastrophe and $n$ have
different parity.
  \label{dd_2}
\end{lem}
\noindent {\it Proof.}  For any Dyck path in $\mathcal{D}''$ having
the catastrophe of size $k$, the number of $U$ minus the number of $D$
equals $k$, which implies that the number of $U$ and $D$ (which is
equal to $n-1$) has the same parity as $k$.\hfill $\Box$

\medskip

  In the following, for a given path $P$,
  a $D$-step will be called {\em isolated}
  if and only if it does not lie 
  in an occurrence of $DD$.

\begin{figure}[H]
  \begin{center}
    \includegraphics[height=6em]{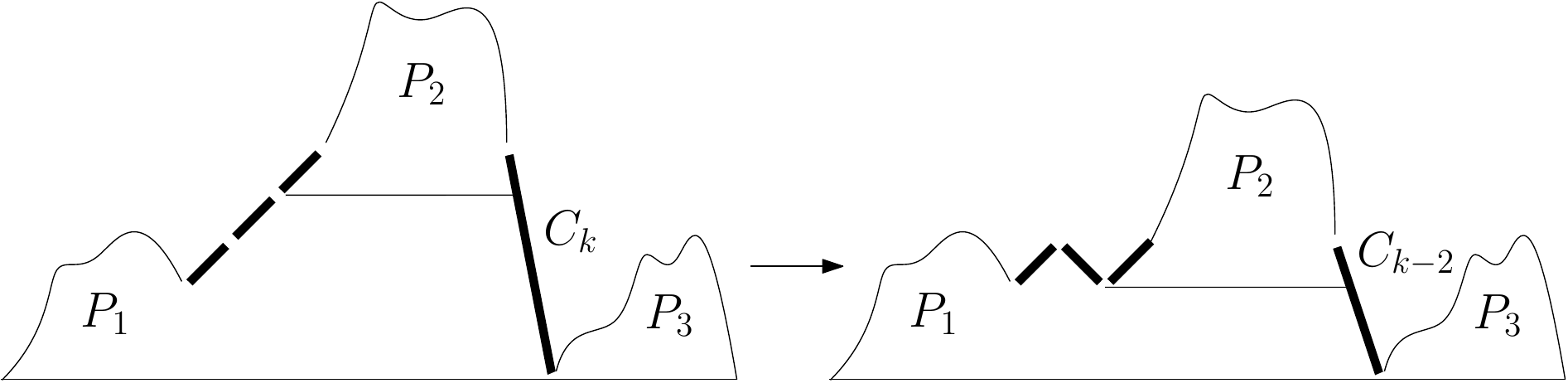}
  \end{center}
  \caption{Illustration for Lemma~\ref{dd_4}.}
  \label{dd_4f}
\end{figure}

\begin{lem}
   Let $P$ be a path in $\mathcal{D}''$ having its catastrophe
  of size $k\geq 2$. We have $P=QC_k R$ where
  $R$ is a Dyck
   path. If $R$ contains isolated $D$-steps, then there is a path in
   $\mathcal{D}''$, $DD$-equivalent to $P$,
   with a catastrophe
    $C_{k+2}$ of size $k+2$, and such that every $D$-step on the right
   of $C_{k+2}$ is not isolated.
 \label{dd_3}
\end{lem}
\noindent {\it Proof.} 
First, we construct $P' = Q' C_k R'$,
where $Q' C_k$ is the result of application 
of Proposition~\ref{prop1} on $QC_k$ and $R'$
is obtained applying $\phi$
from~\cite{Bari1} on $R$. It implies that $P'$ and $P$ 
are $DD$-equivalent and this construction
ensures that the height of the rightmost isolated $D$ step in $P'$  is at most one.
After this, we distinguish three cases according to the position
and the height of the rightmost isolated $D$-step.
\begin{itemize}
  \item[-] if $R'=R_1UD$, then the paths $Q' UR_1UC_{k+2}$ and $P$ are
    equivalent;
    \begin{center}
      \includegraphics[height=6em]{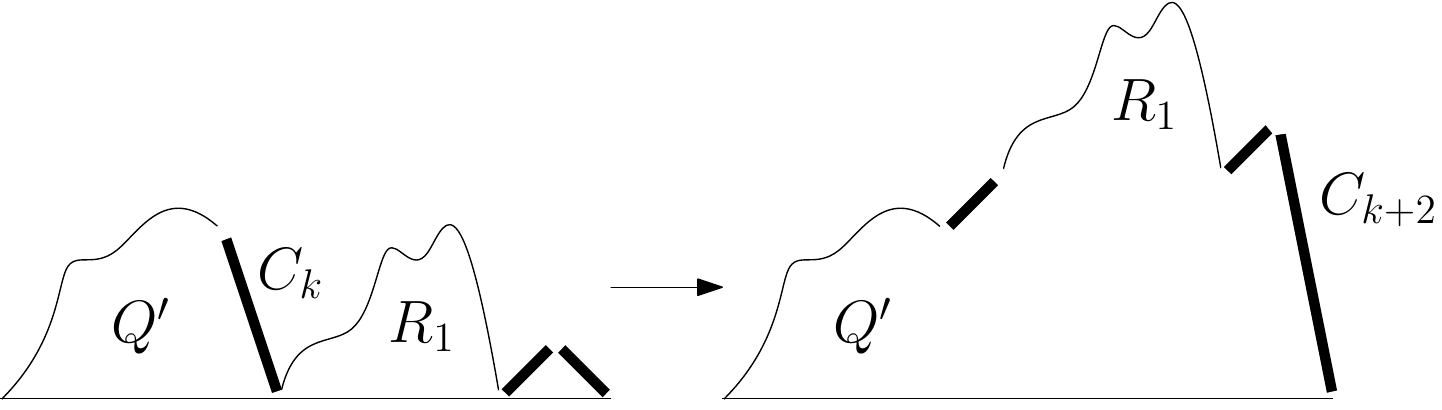}
    \end{center}
  \item[-] if $R' =R_1UDU R_2$ where $UR_2$ is a Dyck path, then $Q' UR_1UC_{k+2}UR_2$ and $P$ are equivalent;
    \begin{center}
      \includegraphics[height=6em]{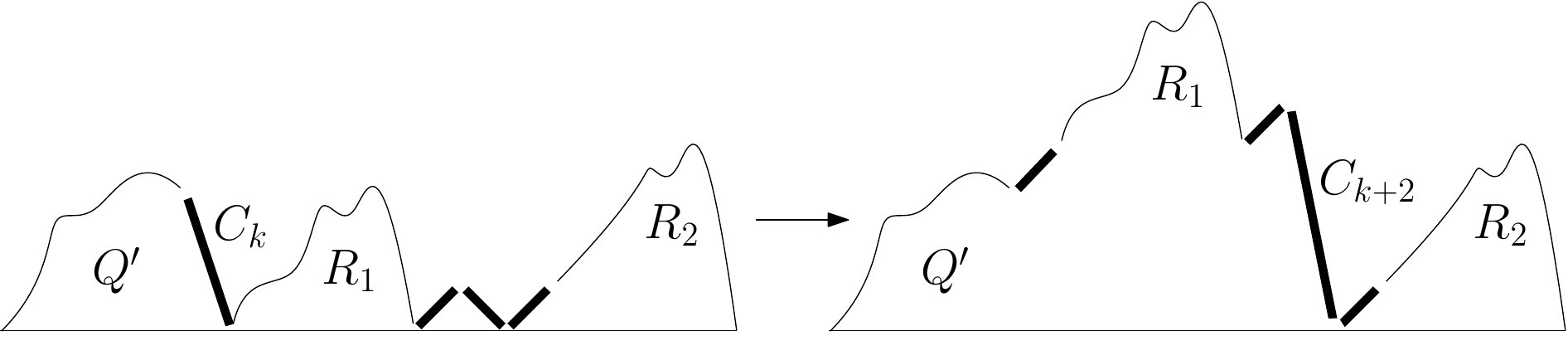}
    \end{center}
  \item[-] if $R' = R_1 UDU R_2$ where $UUR_2$ is a Dyck path then $Q' UR_1C_{k+2}UUR_2$ and $P$ are equivalent.
    \begin{center}
      \includegraphics[height=6em]{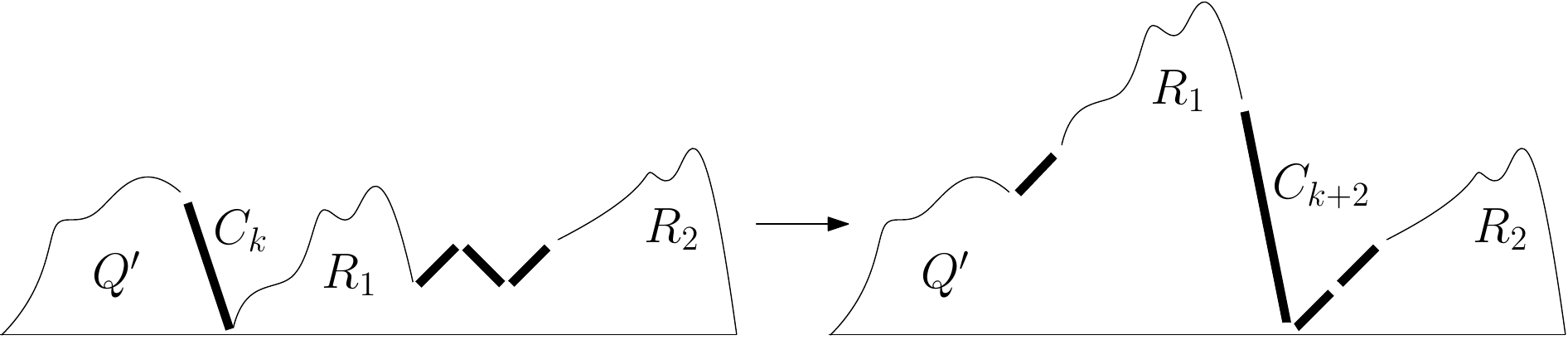}
    \end{center}
\end{itemize}
\hfill $\Box$
\begin{lem} Let $P$ be a path in $\mathcal{D}''$ satisfying
    Condition $(C)$ and having its catastrophe of size $k\geq 4$. Then,
    there is a path $Q\in \mathcal{D}''$, equivalent to $P$,
    having its catastrophe of size $k-2$, and such that $C_k$ and
    $C_{k-2}$ are in the same position in $P$ and $Q$, respectively.
  \label{dd_4}
\end{lem}
\noindent {\it Proof.} Before describing the construction of $Q$, it
is worth noticing the following fact. There is an occurrence of $UUU$
at height $k-3$ on the left of $C_k$ in $P$. Indeed, let us consider
the rightmost $U$ at height $k-3$ on the left of $C_k$; clearly it
necessarily precedes a step $U$ which constitutes an occurrence of
$UU$ at height $k-3$. Moreover, since $P$ avoids $UUDU$, the two steps
$UU$ are necessarily
followed by another step $U$, and $P$ can be
decomposed $P=P_1UUUP_2C_kP_3$ where $P_1$ is a prefix of a Dyck path
ending at height $k-3$, $UP_2D$ is a Dyck path and $P_3$ is a Dyck
path.  We complete the proof by setting $Q=P_1UDUP_2C_{k-2}P_3$. See
Figure~\ref{dd_4f} for an illustration of this construction.\hfill $\Box$

\subsubsection{The case of odd length}

 Let $\mathcal{A}'$ be the subset of paths in $\overline{\mathcal{A}}$ avoiding
 $UDU$, and not ending with $UD$, which also is the set of paths in $\overline{\mathcal{A}}$ without isolated $D$-steps. Let $\mathcal{A}_2$ be the set of paths
 satisfying Condition $(C)$ and having only one catastrophe $C_2$ at
 the end.

\begin{prop}
   A set of representatives of the $DD$-equivalence classes of odd
   length Dyck paths with catastrophes, is the set $\mathcal{R}_1$ of
   paths $R=AA'$ where $A\in\mathcal{A}_2$ and $A'\in\mathcal{A}'$.
   \label{dd_p2}
\end{prop}

\begin{figure}[H]
  \begin{center}
    \includegraphics[height=3em]{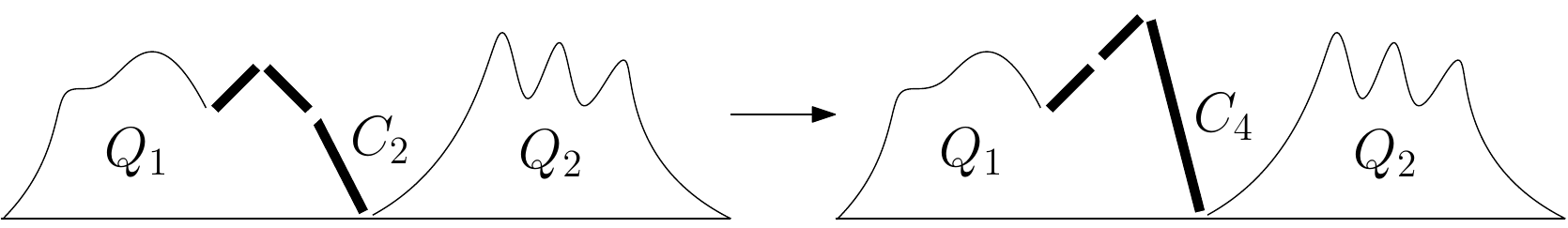}
  \end{center}
  \caption{Illustration for Proposition~\ref{dd_p2}.}
  \label{dd_3f}
\end{figure}

\noindent {\it Proof.} First, let us prove that for any path $P\in
\mathcal{E}$, there is a path $Q\in \mathcal{R}_1$ lying in the same
class.  We obtain $Q$ by applying the following process: ($i$) we
apply Lemma~\ref{dd_1}; ($ii$) we apply Lemma~\ref{dd_3}; ($iii$) we
apply $\psi$ extended to Dyck paths with only one catastrophe (see
Proposition~\ref{prop1}); ($iv$) we apply Lemma~\ref{dd_4} and
Proposition~\ref{prop1}, and we repeat it until the size of the
catastrophe reaches $2$ (thanks to Lemma~\ref{dd_2}). In the case
where $UD$ precedes $C_2$, i.e.~$Q=Q_1UDC_2Q_2$, we replace the path
with $Q_1UUC_4Q_2$ (see Fig.~\ref{dd_3f}) and we apply
Lemma~\ref{dd_4} and Proposition~\ref{prop1}.

 Now, it suffices to prove that two different paths $Q$ and $Q'$ in
 $\mathcal{R}_1$ with the same length cannot belong to the same
 class. Let us assume that $Q$ and $Q'$ lie in the same class. We
 decompose $Q=RC_2S$ and $Q'=R'C_2S'$. Suppose that $|S|\geq
 |S'|$. Since $S'$ does not contain isolated $D$-steps, $S'$ is
 entirely determined by the positions of occurrences of $DD$, which
 implies that $S$ has the form $S=TS'$. Since $S$ does not contain an
 isolated $D$-step, $T$ cannot end with $UD$, and since $Q$ is
 $DD$-equivalent to $Q'$ it cannot end with $DD$, which means that $T$
 is empty. So we have $|S|=|S'|$ and then $S=S'$. Using Proposition~\ref{prop1}
 for $RC_2$ and $R'C_2$, we obtain $Q=Q'$. \hfill $\Box$

 \begin{thm} (The case of odd length) The o.g.f.\ for the number of
   $DD$-equivalence classes of odd length paths in $\mathcal{E}$ is
   given by:

  $$N(x) = \frac{8 x^3}{ (x^2 + \sqrt{1 - 3 x^4 - 2 x^2} + 1)^2 (- x^4 + (1+x^2) \sqrt{1 - 3 x^4 - 2 x^2} - 2 x^2  +
    1 )}.$$

  We have
  $ N(x) =  x^3 + x^5 + 5 x^7 + 11 x^9 + 33 x^{11} + 88 x^{13} + 247 x^{15} + O(x^{16})$.
\label{thm19}
\end{thm}

 \noindent {\it Proof.}
 By Proposition~\ref{dd_p2}, it suffices to provide the o.g.f.\ $N(x)$ for $\mathcal{A}_2 \times \mathcal{A}'$ with respect to the length. Sets  $\mathcal{A}_2$ and $\mathcal{A}'$ are constructed as follows
\begin{center}
  \includegraphics[height=11em]{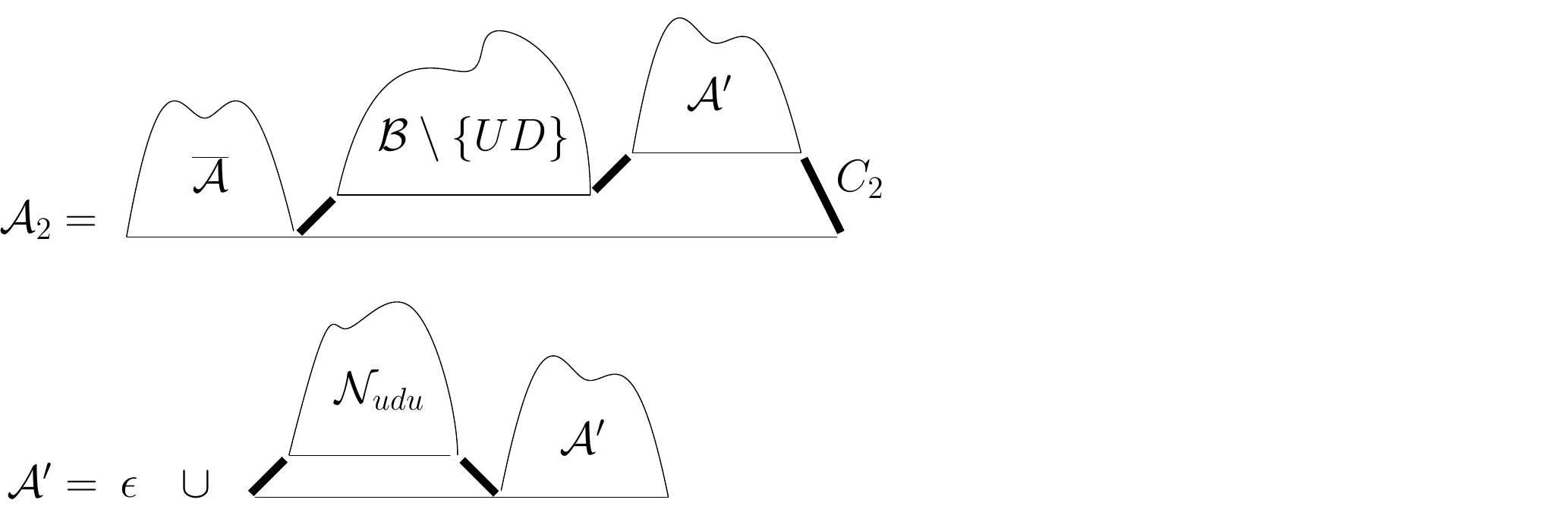}
\end{center}
where $\mathcal{N}_{udu}$ is the set of non-empty Dyck paths avoiding
$UDU$, $\overline{\mathcal{A}}$ is the set of representatives for the $DD$-equivalence
classes of Dyck paths, e.g. Dyck paths avoiding $UUDU$ and having all $UDU$ at height 0 or 1, and $\mathcal{B}$ is the set of Dyck paths having all $UDU$ on $x$-axis
and not starting with $UDU$. From~\cite{Bari1} it follows that

$$
\begin{aligned}
\overline{A}(x) & = \frac{1- x^{2} + \sqrt{- 3 x^{4} - 2 x^{2} + 1} }{1+x^{6} + x^{4} - 3 x^{2} - \left(x^{4} - 1\right) \sqrt{- 3 x^{4} - 2 x^{2} + 1} },\\
B(x) & = \frac{2-x^2- x^{4} + x^{2} \sqrt{- 3 x^{4} - 2 x^{2} + 1}}{1- x^{2} + \sqrt{- 3 x^{4} - 2 x^{2} + 1}},
\end{aligned}
$$
and from~\cite{Bari1,sun} we have $N_{udu}(x)= (x^{2} - \sqrt{- 3 x^{4} - 2 x^{2} + 1} + 1) / (2 x^{2}) - 1$.

 Finally, we compute $N(x) = \overline{A}(x) \left( B(x) - x^2 \right) A'(x) x^3 A'(x)$ where $A'(x)$ satisfies $A'(x)=1+x^2N_{udu}(x)A'(x)$.

 \hfill $\Box$

\subsubsection{The case of even length}

\begin{prop}
  A set of representatives of the $DD$-equivalence classes of even
  length Dyck paths with catastrophes is the set $\mathcal{R}_2$
  defined by the union of $\overline{\mathcal{A}}$ with the sets $\mathcal{S}_0$
  and $\mathcal{S}_1$ defined below:

- $\mathcal{S}_0$ is the set of Dyck paths $S=PDUUC_3$ where $PDU$ is
a prefix of Dyck path satisfying condition $(C)$, and ending at height $2$,

- $\mathcal{S}_1$ is the set of Dyck paths $S=PDDC_3$ where $P$ is a
prefix of Dyck path satisfying condition $(C)$, and ending at height $5$.
 \label{dd_p3}
\end{prop}

\noindent {\it Proof.} Let $P$ be a path in $\mathcal{E}$. If $P$ does
not contain any catastrophe, then using Lemma~1 from~\cite{Bari1},
there is a unique path in $\overline{\mathcal{A}}$ equivalent to
$P$. Now, let us assume that $P$ contains at least one catastrophe. We
obtain $Q$ by applying the following process: ($i$) we apply
Lemma~\ref{dd_1}; we apply Lemma~\ref{dd_3}; we apply $\psi$ (extended
to Dyck paths with only one catastrophe); we apply Lemma~\ref{dd_4} followed by
Proposition~\ref{prop1} until the path contains only one catastrophe $C_3$
(thanks to Lemma~\ref{dd_2}). We decompose the obtained path $Q=RC_3S$ where
$S$ is a Dyck path with no isolated $D$-step. With the same argument
used at the beginning of the proof of Lemma~\ref{dd_4}, $R$ is of the form
$R=R'UUUR''$ where $R'$ and $UR''D$ are Dyck paths.
\begin{enumerate}
    \item [($a$)] if $R''$ is empty, i.e.~$R=R'UUU$, we replace $Q$
with $R'UDUDS\in\overline{\mathcal{A}}$;
\begin{center}
  \includegraphics[height=2.5em]{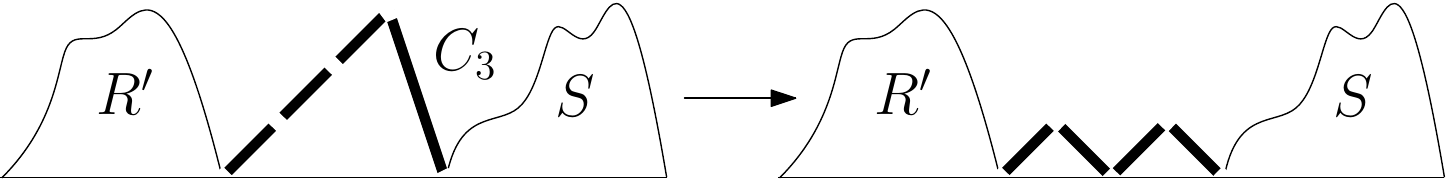}
\end{center}

  \item [($b$)]  if $R''=TU$ ends with $U$, i.e.~$R=R'UUUTU$, we
replace $Q$ with $R'UDUTUDS\in\overline{\mathcal{A}}$;
\begin{center}
  \includegraphics[height=5em]{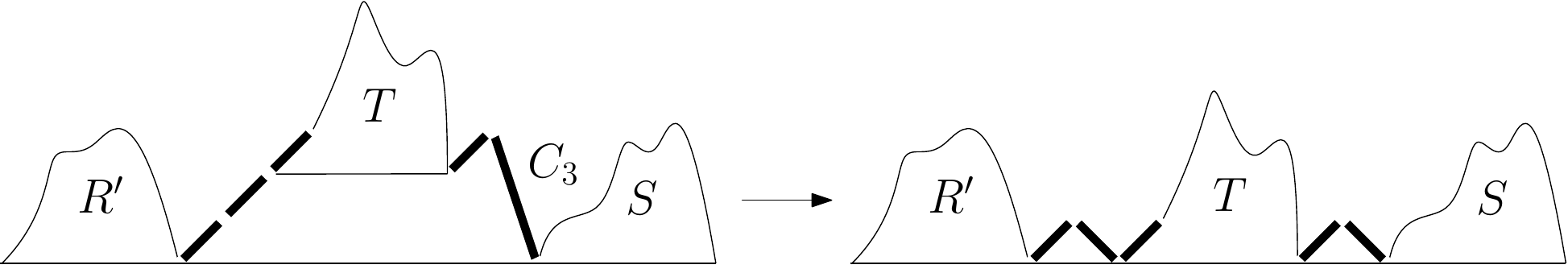}
\end{center}

  \item [($c$)]  if $R''=TUD$ ends with $UD$, i.e.~$R=R'UUUTUD$, we
replace $Q$ with $R'UDUTUUC_3S$ (see the illustration below); In the case where $T$ ends with $U$, we apply the transformation $(a)$ and we obtain a path in $\overline{\mathcal{A}}$;

\begin{center}
  \includegraphics[height=5em]{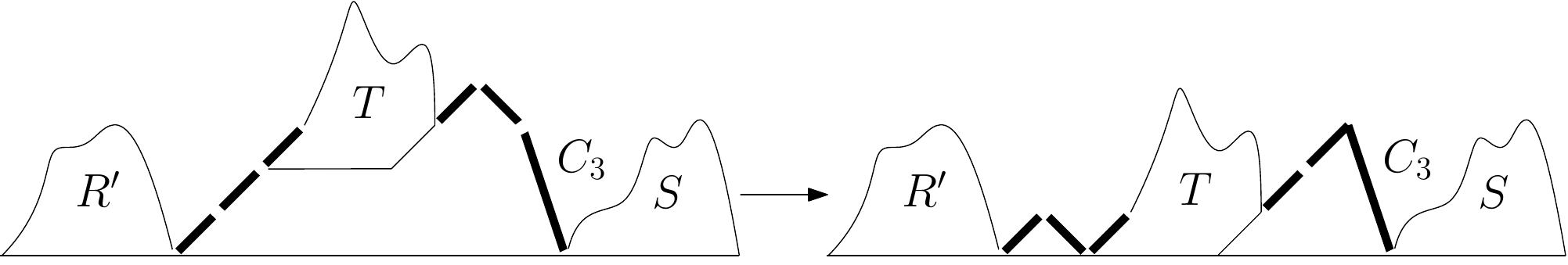}
\end{center}
\end{enumerate}

After applying all these transformations, the path is either in $\overline{\mathcal{A}}$, or it satisfies Condition ($C$) and contains only one catastrophe $C_3$ that follows an occurrence $DUU$ or $DD$.

If $S$ is empty then the path lies in $\mathcal{R}_2$. Otherwise, $S$ does not contain isolated $D$-step and it starts necessarily with $UU$. In this case we set
$S=UUS'$ and make the catastrophe disappear as follows:

\begin{enumerate}
   \item [($d$)]  if $R$ ends at height 3 with $DUU$, we set $R=R'DUU$ and replace $Q$
with $R'DUDUDUS'\in \overline{\mathcal{A}}$;
\begin{center}
  \includegraphics[height=3em]{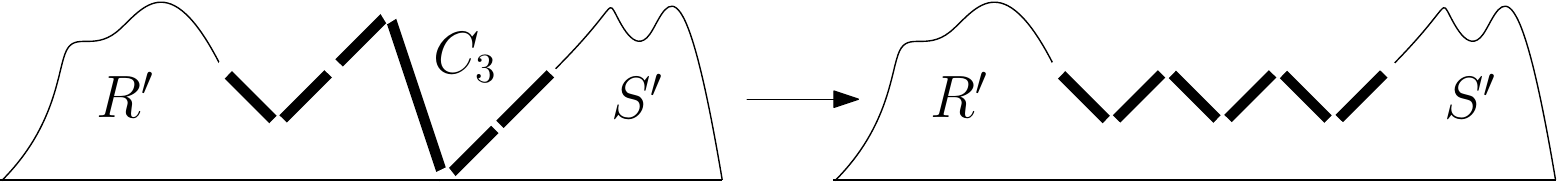}
\end{center}

 \item [($e$)] if $R$ ends with $UUUT$, where $UT$ is a prefix of a Dyck path ending at height 1 and with $DD$,
 we set $R = R'UUUT$ and replace $Q$ by $R'UDUTUDUS'$.
\begin{center}
  \includegraphics[height=6em]{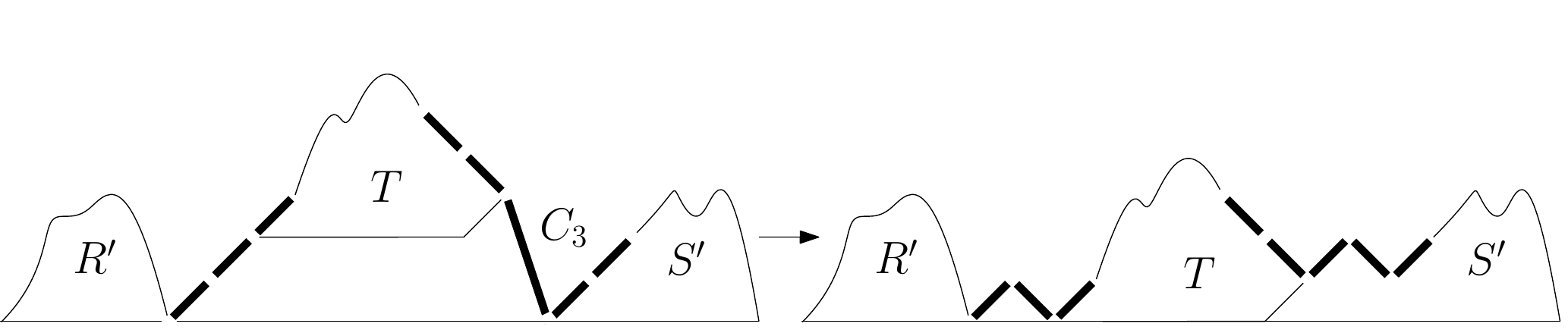}
\end{center}
\hfill $\Box$
 \end{enumerate}

\begin{thm} (The case of even length)  The o.g.f.\ $V(x)$ for the number of $DD$-equivalence classes of even length paths in $\mathcal{E}$ is given by:
$$V(x)=\frac{4\left(x^{4}-x^{2}-1\right) \sqrt{-3 x^{4}-2 x^{2}+1}+8
    x^{8}+20 x^{6}+8 x^{4}-4}{\left(\left(-x^{2}-1\right) \sqrt{-3
      x^{4}-2 x^{2}+1}+x^{4}+2 x^{2}-1\right) \left(x^{2}+\sqrt{-3
      x^{4}-2 x^{2}+1}+1\right)^{2}}.$$

We have $V(x)=1+x^{2}+2 x^{4}+4 x^{6}+11 x^{8}+27 x^{10}+73 x^{12}+194 x^{14}+529 x^{16}+1448 x^{18}+O\! \left(x^{19}\right)$.
\label{thm20}
\end{thm}
 \noindent {\it Proof.}  By Proposition~\ref{dd_p3}, it suffices to
 provide the o.g.f.\ $M(x)$ for
 $\overline{\mathcal{A}}\cup\mathcal{S}_0 \cup \mathcal{S}_1$ with
 respect to the length. Sets $\mathcal{S}_0$ and $\mathcal{S}_1$ are
 constructed as follows:
\begin{center}
  \includegraphics[height=15em]{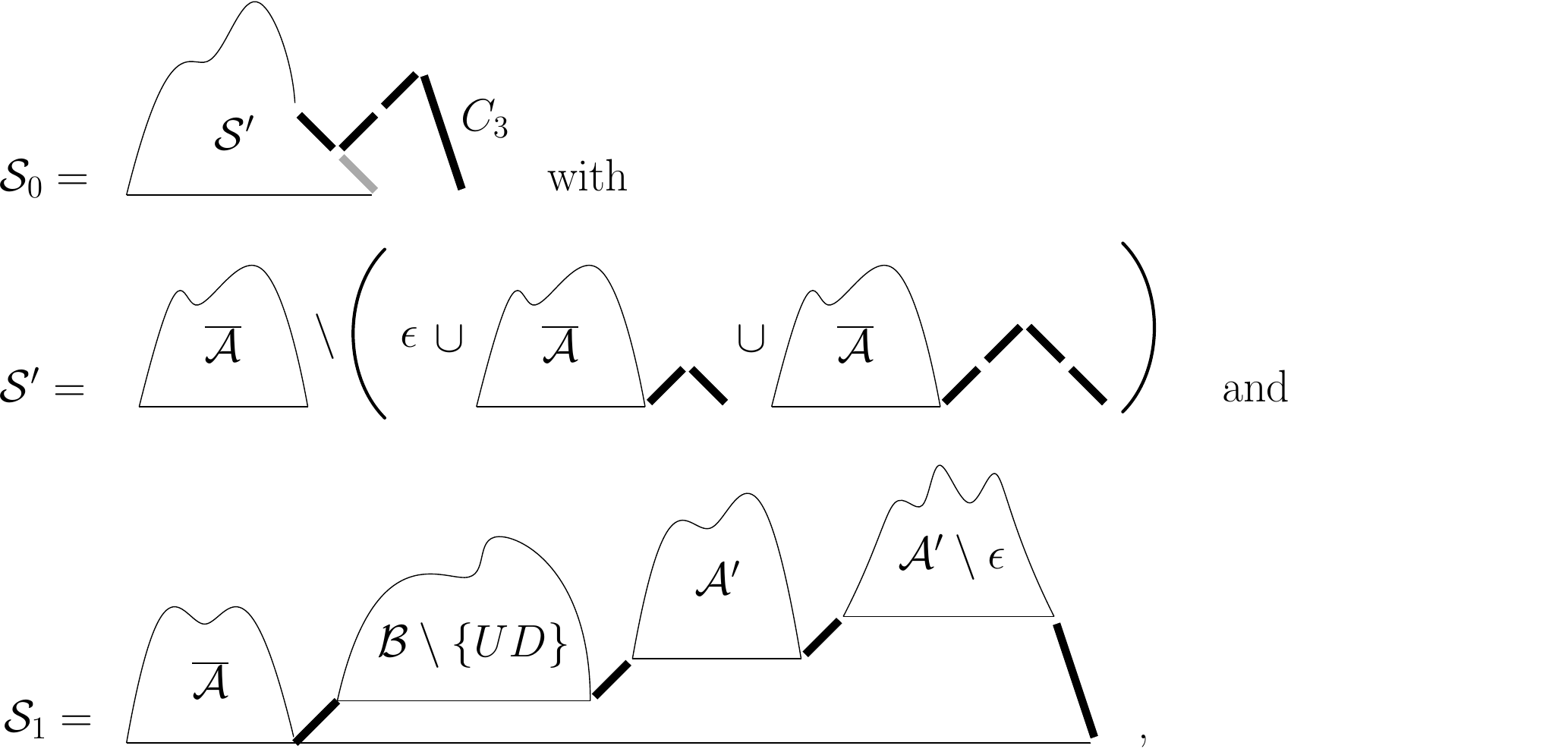}
\end{center}
where $\mathcal{A}'$ is the set of non-empty Dyck paths avoiding
$UDU$ not ending with $UD$, $\overline{\mathcal{A}}$ is the set of representatives for the
$DD$-equivalence classes of Dyck paths, e.g. Dyck paths avoiding
$UUDU$ and having all $UDU$ at height 0 or 1, and $\mathcal{B}$ is the
set of Dyck paths having all $UDU$ on $x$-axis and not starting with
$UDU$.

Then, we
have $$S_0(x)=x^2\cdot(\overline{A}(x)-1-x^2\cdot\overline{A}(x)-x^4\cdot\overline{A}(x))$$

and $$S_1(x)=x^4\cdot\overline{A}(x)(B(x)-x^2)A'(x)(A'(x)-1)$$ which
gives the claimed result after a simple calculation.
\hfill $\Box$

\subsubsection{The general case}
Combining Theorems~\ref{thm19} and~\ref{thm20}, we deduce:

\begin{thm} The o.g.f.\ $L(x)$ for  the number of $DD$-equivalence in $\mathcal{E}$ is given by
$$\frac{\left(4 x^{4}-4 x^{2}-4\right) \sqrt{-3 x^{4}-2 x^{2}+1}+8 x^{8}+20 x^{6}+8 x^{4}-8 x^{3}-4}{\left(\left(-x^{2}-1\right) \sqrt{-3 x^{4}-2 x^{2}+1}+x^{4}+2 x^{2}-1\right) \left(x^{2}+\sqrt{-3 x^{4}-2 x^{2}+1}+1\right)^{2}},
$$
and it satisfies the following equation:
$x^{10}-x^{9}+4 x^{8}-3 x^{7}+5 x^{6}-5 x^{5}+2 x^{4}-2 x^{3}+x^{2}-x +1+\left(-x^{10}+x^{9}-6 x^{8}+5 x^{7}-9 x^{6}+6 x^{5}-3 x^{4}+3 x^{3}+x -1\right) L \! \left(x \right)+\left(x^{10}-x^{9}+4 x^{8}-3 x^{7}+5 x^{6}-3 x^{5}+2 x^{4}-x^{3}\right)L \! \left(x \right)^{2}=0$.

The first coefficients of $x^n$, $n\geq 2$, of the series expansion of $L(x)$  are $1,1,2,1,4,5,11,11,27,33,73,88,194, 247,\ldots$, and they do not correspond to a part of a sequence listed in~\cite{OEIS}.

An asymptotic approximation of these coefficients is
$$\frac{41\sqrt{2}}{2\sqrt{\pi}}\cdot
\frac{\sqrt{3}^{n+1}}{n^{\frac{3}{2}}} \mbox{ for } n \mbox{ even, and } $$
$$\frac{135\sqrt{2}}{4\sqrt{\pi}} \cdot
\frac{\sqrt{3}^n}{n^{\frac{3}{2}}} \mbox{ otherwise. }
$$
\label{thm21}
\end{thm}

  \paragraph{Acknowledgments.}
The authors would like to thank the anonymous referees for useful
remarks and comments.  This work was supported in part by the project
ANER ARTICO funded by Bourgogne-Franche-Comté region of France and ANR
JCJC PICS.


\begin{thebibliography}{99}
\bibitem{Ban} C. Banderier and M.~Wallner, Lattice paths with catastrophes, \emph{Electron. Notes Discret. Math.}, \textbf{59} (2017), 131--146.

\bibitem{Bari1} J.-L.~Baril and A.~Petrossian,  Equivalence of Dyck paths modulo some statistics, \emph{Discrete Math.}, \textbf{338} (2015), 655--660.

\bibitem{Bari2} J.-L.~Baril and  A.~Petrossian, Equivalence classes  of Motzkin paths modulo a pattern of length at most two, \emph{J. Integer Seq.}, \textbf{18} (2015), 15.7.1.

\bibitem{Bari4} J.-L. Baril and S.~Kirgizov,  Bijections from Dyck and Motzkin meanders with catastrophes to pattern avoiding Dyck paths,  \emph{Discrete Math. Lett.}, \textbf{7} (2021), 5--10.

\bibitem{Bari3} J.-L.~Baril, S.~Kirgizov, and A.~Petrossian, Enumeration of \L{}ukasiewicz paths modulo some patterns, \emph{Discrete Math.}, \textbf{342}(4) (2019), 997--1005.


\bibitem{Barasii} J.-L.~Baril, J.L.~Ram\'{\i}rez, L.M.~Simbaqueba,  Equivalence classes of skew Dyck paths Modulo some patterns, \emph{Integers}, \textbf{22} (2022).

\bibitem{Deut} E.~Deutsch,  Dyck path enumeration, \emph{Discrete Math.}, \textbf{204}(1999), 167--202.


\bibitem{flajolet} P.~Flajolet and R.~Sedgewick,  \emph{Analytic Combinatorics}, Cambridge University Press,  2009.



\bibitem{Kri} A. Krinik, G. Rubino, D. Marcus, R.J. Swift, H. Kasfy,  H. Lam.
\newblock Dual processes to solve single server systems.
\newblock {\em J. Stat. Plan. Inference}, 135(2005), 1, 121--147.


\bibitem{Manes}  K.~Manes, A.~Sapounakis, I.~Tasoulas, and P.~Tsikouras, Equivalence classes of ballot paths modulo strings of length 2 and 3, \emph{Discrete Math.}, \textbf{339}(10) (2016), 2557--2572.

\bibitem{orlov} A.G.~Orlov, On asymptotic behavior of the Taylor coefficients of algebraic functions, \emph{Sib. Math. J.}, \textbf{25}(5) (1994), 1002--1013.

\bibitem{OEIS} N.~J.~A.~Sloane, The On-Line Encyclopedia of Integer Sequences. Available at \url{https://oeis.org/}.

\bibitem{sun} Y.~Sun, The statistic {``}number of udu’s{"} in Dyck paths,  \emph{Discrete Math.}, \textbf{287} (2004), 177--186.

\end{thebibliography}
\end{document}